\documentclass[twoside]{article}

\usepackage{microtype}

\usepackage[accepted]{aistats2021}
\usepackage[round]{natbib}

\usepackage{float}
\usepackage[dvipsnames]{xcolor}
\usepackage{amsmath,amsfonts,amssymb,amsthm}
\usepackage{mathtools}
\usepackage{enumerate}
\usepackage[capitalise]{cleveref}
\usepackage{lipsum}
\usepackage{caption}
\usepackage{subcaption}
\usepackage{float}
\usepackage{url}

\usepackage{tikz}
\usepackage{pgfplots}
\usetikzlibrary{backgrounds}
\usetikzlibrary{intersections}
\usepgfplotslibrary{fillbetween}

\usetikzlibrary{external}   
\tikzsetexternalprefix{external/}  

\newlength{\figheight}
\newlength{\figwidth}

\newcommand{\T}{^{\intercal}}
\renewcommand{\d}{\:\mathrm{d}}
\DeclareMathOperator*{\argmax}{arg\,max}

\DeclareMathOperator{\Cov}{Cov}
\DeclareMathOperator{\diag}{diag}

\newtheorem{proposition}{Proposition}
\newtheorem{remark}{Remark}
\renewcommand{\Sigma}{\varSigma}

\definecolor{notecolor}{RGB}{0,82,147}

\begin{document}

%

%
\runningauthor{Nathanael Bosch, Philipp Hennig, Filip Tronarp}

\twocolumn[

\aistatstitle{Calibrated Adaptive Probabilistic ODE Solvers}

\aistatsauthor{ Nathanael Bosch\textsuperscript{1} \And Philipp Hennig\textsuperscript{1,2} \And  Filip Tronarp\textsuperscript{1} }

\aistatsaddress{
  \And
  \textsuperscript{1}University of Tübingen\\
  \textsuperscript{2}Max Planck Institute for Intelligent Systems, Tübingen, Germany
  \\ \texttt{\{nathanael.bosch, philipp.hennig, filip.tronarp\}@uni-tuebingen.de}
  \And
} ]

\begin{abstract}
  Probabilistic solvers for ordinary differential equations assign a posterior measure to the solution of an initial value problem.
  The joint covariance of this distribution provides an estimate of the (global) approximation error.
  The contraction rate of this error estimate as a function of the solver’s step size identifies it as a well-calibrated worst-case error,
  but its explicit numerical value for a certain step size
  is not automatically a good estimate of the explicit error.
  Addressing this issue, we introduce, discuss, and assess several probabilistically motivated ways to calibrate the uncertainty estimate.
  Numerical experiments demonstrate that these calibration methods interact efficiently with adaptive step-size selection, resulting in descriptive, and efficiently computable posteriors.
  We demonstrate the efficiency of the methodology by benchmarking against the classic, widely used Dormand--Prince 4/5 Runge--Kutta method.
\end{abstract}

\section{INTRODUCTION}

\begin{figure*}[t]
  \centering
  \begin{subfigure}[b]{0.48\linewidth}
    \centering
    \includegraphics[width=\linewidth]{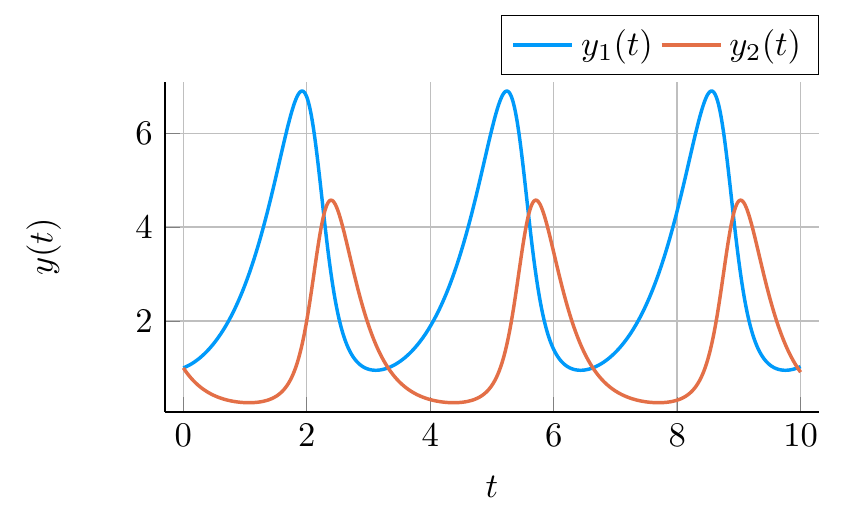}%
  \end{subfigure}
  ~
  \begin{subfigure}[b]{0.48\linewidth}
    \centering
    \includegraphics[width=\linewidth]{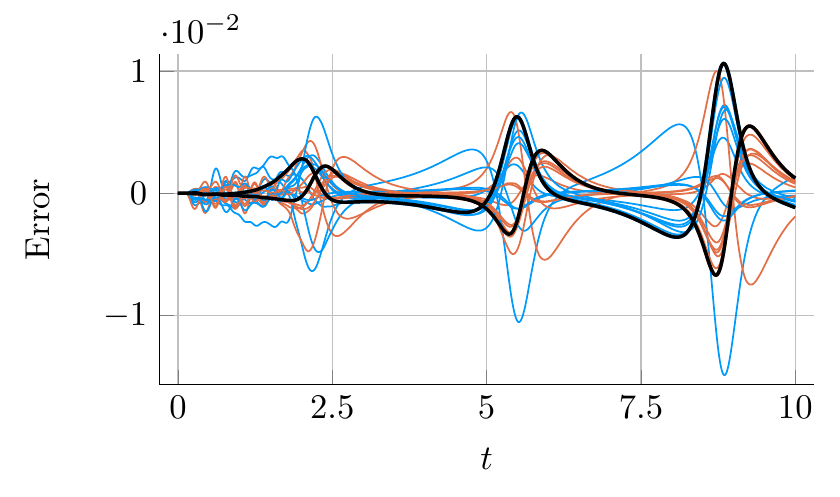}%
  \end{subfigure}
  \caption{
    \emph{Probabilistic solution of the Lotka-Volterra equations}.
    Left:
    Posterior mean returned by the probabilistic solver.
    For the chosen tolerance levels the returned uncertainties are too small to be visually separated from the mean.
    Right:
    True error trajectory (black) and sampled error trajectories (colored).
    Both the true solution and the samples exhibit similar patterns, indicating a well-structured and calibrated uncertainty estimate of the provided posterior distribution.
  }
\end{figure*}

Ordinary differential equations (ODEs) arise in almost all areas of science and engineering.
In the field of machine learning,
recent work on normalizing flows
\citep{pmlr-v37-rezende15}
and neural ODEs
\citep{DBLP:conf/nips/ChenRBD18}
lead to a particular surge of interest.
In this paper we consider
initial value problems (IVPs),
defined by an ODE
\begin{equation}
  \label{eq:ivp}
  \dot{y}(t) = f \left( y(t), t \right), \qquad \forall t \in [t_0, T],
\end{equation}
with vector field \(f: \mathbb{R}^d \times \mathbb{R} \to \mathbb{R}^d\) and initial value \(y(t_0) = y_0 \in \mathbb{R}^d\).

Limited by finite computational resources, the numerical solution of an IVP is inevitably only an approximation.
Though, most classic numerical solvers do not return an estimate of their own numerical error, leaving it to the practitioner to evaluate the reliability of the result.
The field of \emph{probabilistic numerics} (PN)
\citep{HenOsbGirRSPA2015,Oates_2019}
seeks to overcome this ignorance of numerical uncertainty.
By treating numerical algorithms as problems of statistical inference, the numerical error can be quantified probabilistically.

One particular class of probabilstic numerical solvers for ODEs treats IVPs as a Gauss--Markov regression problem
\citep{DBLP:journals/sac/TronarpKSH19,tronarp20_bayes_ode_solver},
the solution of which can be efficiently approximated with Bayesian filtering and smoothing
\citep{DBLP:books/daglib/0039111}.
These so-called ODE filters relate to classic multistep methods
\citep{schober16_probab_model_numer_solut_initial_value_probl},
and have been shown to converge to the true solution of the IVP with high polynomial rates while providing (asymptotically) well-calibrated confidence intervals
\citep{Kersting2020}.
In practice, however, there remain two gaps:
First, efficient implementation of ODE solvers require adaptive step-size selection, which has not received much attention in the past;
second, the calibration of the posterior uncertainty estimates depends on the choice of specific diffusion hyperparameters.
Both are addressed in this paper.

The contributions of this paper are the following:
We introduce and discuss uncertainty calibration methods for models with both constant and time-varying diffusion, and extend existing approaches with multivariate parameter estimates.
After calibration, the probabilistic observation model provides an objective for local error control and adaptive step-size selection, enabling the solvers to make efficient use of their computational budget.
The resulting probabilistic numerical ODE solvers are evaluated and compared for a large range of configurations and tolerance levels,
demonstrating descriptive posteriors and computational efficiency comparable to the classic Dormand--Prince 4/5 Runge--Kutta method.

\section{PROBABILISTIC ODE SOLVERS}
\label{sec:solvers}
In this paper, the solution of the IVP is posed as a Bayesian inference problem
\citep{cockayne17_bayes_probab_numer_method}.
By considering Gauss--Markov priors on the state-space, the problem is reduced to a case of Bayesian state estimation
\citep{DBLP:journals/sac/TronarpKSH19}.
In the following, we introduce the state estimation problem and describe the approximate inference procedure.

\subsection{ODE Solutions as State Estimation}
\label{sec:continuous-time-model}
\emph{A priori}, we model the solution together with \(q \in \mathbb{N}\) of its derivatives
as a Gauss--Markov process
\({X(t) = \left[\left(X^{(0)}(t)\right)\T, \left(X^{(1)}(t)\right)\T, \dots, \left(X^{(q)}(t)\right)\T\right]\T}\),
where \(X^{(i)}(t)\) models the \(i\)-th derivative \(y^{(i)}(t)\).
More precisely, \(X(t)\) is the \(q\)-times integrated Wiener process (IWP), which solves the stochastic differential equations
\begin{subequations}
  \begin{align}
    \d X^{(i)}(t) &= X^{(i+1)}(t) \d t, \quad i=1, \dots, q-1\\
    \d X^{(q)}(t) &= \Gamma^{1/2} \d B(t), \\
    X(t_0) &\sim \mathcal{N}(\mu_0, \Sigma_0),
  \end{align}
  \label{eq:prior}
\end{subequations}
where \(\Gamma^{1/2}\) is the symmetric square root of some positive semi-definite matrix \(\Gamma \in \mathbb{R}^{d \times d}\).

The continuous-time model can be described by transition densities
\citep[Section 6.2]{sarkka_solin_2019}
\begin{equation}
  X(t + h) \,|\, X(t) \sim \mathcal{N} \left( A(h) X(t), Q(h) \right).
  \label{eq:discrete-prior}
\end{equation}
For the IWP prior,
\({A(h) \in \mathbb{R}^{d(q+1) \times d(q+1)}}\) and \({Q(h) \in \mathbb{R}^{d(q+1) \times d(q+1)}}\)
are of the form
\begin{subequations}
  \begin{align}
    A(h) &= \breve{A}(h) \otimes I_d, \\
    Q(h) &= \breve{Q}(h) \otimes \Gamma,
  \end{align}
  \label{eq:AQ-kron}
\end{subequations}
with \(\breve{A}(h), \breve{Q}(h)\) given by \citet[Appendix A]{Kersting2020}
\begin{subequations}
  \begin{align}
    \breve{A}_{ij}(h) &= \mathbb{I}_{i \leq j} \frac{h^{j-1}}{(j-i)!}, \\
    \breve{Q}_{ij}(h) &= \frac{h^{2q+1-i-j}}{(2q+1-i-j)(q-i)!(q-j)!}.
  \end{align}
  \label{eq:A0Q0}
\end{subequations}

To relate the prior to the solution of the IVP, define the measurement process
\begin{equation}
Z(t) = X^{(1)}(t) - f \left( X^{(0)}(t) \right).
\end{equation}
The probabilistic numerical solution of the ODE is computed by conditioning \(X(t)\) on the event that the realisation $z(t)$ of $Z(t)$ is zero on the grid $\{t_n\}_{n=1}^N$
\citep{DBLP:journals/sac/TronarpKSH19}
\begin{equation}
  z_n := z(t_n) = 0, \quad n = 1,\ldots,N.
\end{equation}
The resulting inference problem of computing
\begin{equation}
  \label{eq:posterior}
  p \left( X(t) \,|\, \{z_n\}_{n=1}^N \right)
\end{equation}
is, for non-linear vector fields \(f\), known as a non-linear Gauss--Markov regression problem and is in general intractable,
but it is possible to efficiently compute approximations \citep{DBLP:books/daglib/0039111}.

\subsection{Approximate Gaussian Inference}
We consider approximate Bayesian inference based on linearization by Taylor-series expansion, known as the extended Kalman filter (EKF) in statistical signal processing
\citep[Section 5.2]{DBLP:books/daglib/0039111}.
By linearizing the measurement likelihood we can efficiently and iteratively compute approximations
\begin{subequations}
\begin{align}
  p \left( X(t_n) \,|\, \{z_i\}_{i=1}^{n-1} \right) &\approx \mathcal{N} \left( \mu_n^P, \Sigma_n^P \right), \\
  p \left( X(t_n) \,|\, \{z_i\}_{i=1}^n \right) &\approx \mathcal{N} \left( \mu_n^F, \Sigma_n^F \right), \\
  p \left( z_n \,|\, \{z_i\}_{i=1}^{n-1} \right) &\approx \mathcal{N} \left( \hat{z}_n, S_n \right),
\end{align}
\end{subequations}
through the following \emph{prediction} and \emph{update} steps
\citep[Section 5.2]{DBLP:books/daglib/0039111}:

Prediction:
  \begin{subequations}
    \begin{align}
      \mu_n^P &= A(h_{n-1}) \mu_{n-1}^F, \\
      \Sigma_n^P &= A(h_{n-1}) \Sigma_{n-1}^F A(h_{n-1})\T + Q(h_{n-1}).
    \end{align}
  \end{subequations}

Update:
\begin{subequations}
  \begin{align}
    \hat{z}_n &= E_1 \mu_n^P - f\left( E_0 \mu_n^P, t_n \right), \\
    S_n &= H_n \Sigma_n^P H_n\T, \\
    K_n &= \Sigma_n^P H_n\T S_n^{-1}, \\
    \mu_n^F &= \mu_n^P + K_n (z_n - \hat{z}_n), \\
    \Sigma_n^F &= \Sigma_n^P - K_n S_n K_n\T,
  \end{align}
  \label{eq:update}
\end{subequations}
where \(H_n\) can be either
\( H_n := E_1 \) for a zeroth order approximation, or
\( H_n := E_1 - J_f(E_0 \mu_n, t_n) E_0 \) for a first order approximation of the vector field \(f\),
with \(E_0 := e_0\T \otimes I \), \(E_1 := e_1\T \otimes I \).

The zeroth and first order linearizations correspond to the updates by
\citet{schober16_probab_model_numer_solut_initial_value_probl}
and
\citet{DBLP:journals/sac/TronarpKSH19}, respectively.
In the sequel, we refer to the algorithm with zeroth and first order linearization as EKF0 and EKF1, respectively.

\begin{remark}
  While most classic ODE solvers do not use the Jacobians of the vector field \(f\),
  they play a central role in Rosenbrock methods
  \citep{DBLP:journals/cj/Rosenbrock63, DBLP:journals/siamnum/HochbruckOS08},
  a class of semi-implicit solvers for stiff ODEs
  \citep[Chapter IV.7]{hairer1987solving}.
  In probabilistic solvers, the Jacobian was used in
  a probabilistic multistep method
  \citep{DBLP:conf/nips/TeymurZC16}
  and, more recently, with extended Kalman filtering and smoothing
  \citep{DBLP:journals/sac/TronarpKSH19,tronarp20_bayes_ode_solver}.
\end{remark}

The Bayesian \emph{filtering} posterior for \(X(t)\)
is conditioned only on the measurements obtained before and at the time step \(t\), but does not include future measurements.
Computing the (approximate) full marginal posterior
\(p \left( X(t_n) \,|\, \{z_n\}_{n=1}^{N} \right)\)
can be done with Bayesian \emph{smoothing}.
The extended Rauch--Tung--Striebel smoother,
also called the extended Kalman smoother (EKS),
describes an algorithm to efficiently compute Gaussian approximations
\begin{align}
  p \left( X(t_n) \,|\, \{z_i\}_{i=1}^{N} \right) &\approx \mathcal{N} \left( \mu_n^S, \Sigma_n^S \right)
\end{align}
with a backwards recursion, given by the
\emph{smoothing} step
\begin{subequations}
\begin{align}
  G_n &= \Sigma_n^F A(h_n) (\Sigma_{n+1}^P)^{-1}, \\
  \mu_n^S &= \mu_n^F + G_n (\mu_{n+1}^S - \mu_{n+1}^P), \\
  \Sigma_n^S &= \Sigma_n^F + G_n (\Sigma_{n+1}^S - \Sigma_{n+1}^P) G_n\T.
\end{align}
\end{subequations}
See also \citet[Chapter 9]{DBLP:books/daglib/0039111}.
The approximate posterior for off-the-grid time steps \(t \in [t_0, T]\),
relating to \emph{dense output} in classic numerical solvers
\citep[Chapter II.6]{hairer2008solving},
can be straight-forwardly computed by using the Gauss--Markov property.
In a similar manner as above, we refer to the resulting algorithms with linearization of order zero and one as EKS0 and EKS1, respectively.

The question remains how to set the initial state
\(X_0 \sim \mathcal{N}\left( \mu_0, \Sigma_0 \right)\).
This problem can also be observed in classic multistep methods:
In addition to the multistep formula, they specify a starting procedure using, for example,
Taylor series expansion
\citep{bashforth1883attempt}, one-step methods
\citep[Chapter III.1]{hairer2008solving},
or other iterative procedures \citep{Nordsieck1962OnNI}.
In a similar effort,
\citet{schober16_probab_model_numer_solut_initial_value_probl}
describe an initialization for a probabilistic ODE solver via Runge-Kutta methods.
Since the probabilistic formulation enables us to explicitly quantify uncertainty over the initial values, it is also possible to set the initial state to a zero mean and unit variance Gaussian distribution
and to condition on the correct initial values
\(X^{(0)}(t_0) = y_0\) and \(X^{(1)}(t_0) = f(y_0, t_0)\)
\citep{DBLP:journals/sac/TronarpKSH19, tronarp20_bayes_ode_solver}.
For our experiments, where orders \(q \leq 5\) are used, we found it feasible to compute all derivatives of the correct initial value via automatic differentiation, and we explicitly set
\(\mu_0 = (y(t_0), \dot{y}(t_0), ..., y^{(q)}(t_0))\) and zero covariance.

\begin{remark}
  The exact initial derivatives can be computed efficiently with Tylor-mode automatic differentiation
  \citep{Griewank2000EvaluatingD,bettencourt2019taylor}.
  A more extensive description of an initialization procedure, together with further considerations for a numerically stable implementation, is provided by
  \citet{kraemer2020stable}.
\end{remark}

\section{UNCERTAINTY CALIBRATION}
\label{sec:uncertaintycalibration}
The probabilistic solver with IWP prior presented in
\cref{sec:continuous-time-model}
contains a free parameter $\Gamma$,
which
is of particular importance for the posterior uncertainty, as
it determines the gain of the Wiener process entering the system
in \cref{eq:prior}.
In this section, we present different diffusion models and discuss approaches for estimating this parameter and thereby calibrating uncertainties.

\subsection{Time-Fixed Diffusion Model}
A common approach for Bayesian model selection and parameter estimation is to maximize the marginal likelihood, or \emph{evidence}, of the observed data \(z_{1:n}\), given by
the prediction error decomposition
\citep{schweppe1965evaluation}
\begin{equation}
  p(z_{1:n}) = p(z_1) \prod_{i=2}^{n} p(z_{i} \,|\, z_{1:i-1}).
\end{equation}
For affine vector fields, the Kalman filter computes the marginals
\(p(z_{i}\,|\,z_{1:i-1})\)
exactly, but for non-affine vector fields we solve the Bayesian filtering problem only approximately.
Nevertheless, it is a natural choice to approximate the marginal likelihoods in the same way as the filtering solution is approximated, i.e.~with extended Kalman filtering,
as
\begin{subequations}
  \begin{align}
    \label{eq:quasi-mle}
    p(z_{1:n}) &\approx \prod_{i=1}^{n} \mathcal{N}\left(z_{i}; \hat{z}_i, S_i \right), \\
    \hat{z}_i &= E_1 \mu_i^P - f\left( E_0 \mu_i^P, t_i \right), \\
    S_i &= H_i \Sigma_i^P H_i\T.
  \end{align}
\end{subequations}
Maximizing \cref{eq:quasi-mle} is referred to as \emph{quasi maximum likelihood estimation} in signal processing
\citep{lindstroem}.

\label{sec:uc-fixedmle}
For the case of scalar matrices
\(\Gamma=\sigma^2 I_d\),
\citet[Proposition 4]{DBLP:journals/sac/TronarpKSH19}
provide a (quasi) maximum-likelihood estimate (quasi-MLE) of
\(\sigma^2\), denoted by \(\hat{\sigma}_N^2\).
Assuming an initial covariance of the form
\(\Sigma_0 = \sigma^2 \breve{\Sigma_0}\), \(\hat{\sigma}_N^2\) is given by
\begin{equation}
  \label{eq:uc-fixedmle}
  \hat{\sigma}_N^2 = \frac{1}{Nd} \sum_{n=1}^N (z_n - \hat{z}_n)\T S_n^{-1} (z_n - \hat{z}_n).
\end{equation}
This estimation can be performed on-line in order to provide calibrated uncertainty estimates during the solve, which are required for step-size adaptation.

\label{sec:uc:fixedMV}
The IWP prior with scalar diffusion \(\Gamma = \sigma^2 I_d\) describes the same model for each dimension.
Furthermore, as the measurement matrix \(H_n\) associated with the EKF0 does not depend on the vector field, the estimated uncertainties of each dimension will be the same.
To fix this shortcoming of the EKF0 we propose a model with diagonal \(\Gamma=\diag(\sigma_1^2, \dots, \sigma_d^2)\).
Its quasi-MLE is provided in \Cref{prop:fixedMV} below.

\begin{proposition}
  \label{prop:fixedMV}
  Let
  \(\Gamma=\diag(\sigma_1^2, \dots, \sigma_d^2)\)
  and
  \(\Sigma_0 = \breve{\Sigma}_0 \otimes \Gamma\).
  Then the prediction and filtering covariances computed by the EKF0 with IWP prior are of the form
  \(\Sigma_n^P = \breve{\Sigma}_n^P \otimes \Gamma\),
  \(\Sigma_n^F = \breve{\Sigma}_n^F \otimes \Gamma\),
  and the approximated measurement covariances are given by
  \(S_n = \breve{s}_n \cdot \Gamma\), where \(\breve{s}_n\) is \(\breve{s}_n := e_2\T \breve{\Sigma}_n^P e_2 \).
  The quasi maximum-likelihood estimate of \(\Gamma\), denoted by \(\hat{\Gamma}\), is diagonal and given by
  \begin{equation}
    \hat{\Gamma}_{ii} = \frac{1}{N} \sum_{n=1}^N \frac{(\hat{z}_n)_i^2}{\breve{s}_n}, \qquad i \in \{1, \dots, d\}.
  \end{equation}
\end{proposition}
The proof follows the idea of
\citet[Proposition 4]{DBLP:journals/sac/TronarpKSH19}.
A more detailed derivation can be found in
\cref{sec:proofs:ekf0-fixed}.

\subsection{Time-Varying Diffusion Model}
\label{sec:uc-schober}
To allow for greater flexibility,
\citet{schober16_probab_model_numer_solut_initial_value_probl} propose a model in which \(\Gamma = \Gamma_n\) is allowed to vary for different integration steps \(t_n\).
In such a model, all measurements \(\{z_i\}_{i=1}^{n-1}\) taken before time \(t_n\) are independent of the parameter \(\Gamma_n\).
We obtain
\begin{subequations}
  \label{eq:tv-argmax}
  \begin{align}
      \argmax_{\Gamma_n} p \left( z_{1:n} \right) &= \argmax_{\Gamma_n} p \left( z_n \,|\, z_{1:n} \right) \\
      &\approx \argmax_{\Gamma_n} \mathcal{N} \left( z_n; \hat{z}_n, S_n \right),
  \end{align}
\end{subequations}
with \(S_n = H_n \left[ A_n \Sigma_{n-1}^F A_n\T + (\breve{Q}_n \otimes \Gamma_n) \right] H_n\T\).

To approximately estimate \(\Gamma_n\),
\citet{schober16_probab_model_numer_solut_initial_value_probl}
propose an estimation based on ``local'' errors, a common procedure for error control and step-size adaptation in classic numerical methods
\citep[Chapter II.4]{hairer2008solving}.
Assuming an error-free predicted solution \(\mu_{n-1}^F\) at time \(t_{n-1}\), that is, \(\Sigma_{n-1}^F = 0\), yields
\begin{equation}
  S_n = H_n (\breve{Q}_n \otimes \Gamma_n) H_n\T.
\end{equation}

For scalar matrices
\( \Gamma_n =  \sigma_n^2 I_d \)
this implies
\begin{equation}
  S_n = \sigma_n^2 \cdot H_n (\breve{Q}_n \otimes I_d) H_n\T .
\end{equation}
Computing the quasi-MLE by solving \cref{eq:tv-argmax} yields the
parameter estimate by
\citet{schober16_probab_model_numer_solut_initial_value_probl}:
\begin{equation}
  \hat{\sigma}_n^2 = \frac{1}{d} (z_n - \hat{z}_n)\T \left( H_n (\breve{Q}_n \otimes I_d) H_n\T \right)^{-1} (z_n - \hat{z}_n).
\end{equation}

\label{sec:uc-schobermv}
As for the fixed diffusion model, we can improve the expressiveness of the EKF0 by considering a multivariate model with diagonal
\(\Gamma_n = \diag(\sigma_{n1}^2, \dots, \sigma_{nd}^2)\).
With the local error based estimation, and using
\({H_n = e_1 \otimes I_d}\),
we obtain
\begin{equation}
  S_n = (\breve{Q}_n)_{11} \cdot \Gamma_n.
\end{equation}
With a covariance of this form, \cref{eq:tv-argmax} can be solved and we obtain the quasi-MLE
\begin{equation}
  (\hat{\Gamma}_n)_{ii} = (z_n - \hat{z}_n)_i^2 / (\breve{Q}_n)_{11}, \qquad i \in \{1, \dots, d\},
\end{equation}
as parameter estimate for models with the time-varying, diagonal diffusion.

\section{STEP-SIZE ADAPTATION}
\label{sec:adaptation}
While the algorithm described in
\cref{sec:solvers,sec:uncertaintycalibration}
is able to compute calibrated posterior distributions over solutions to any IVP, it still lacks a common tool to make efficient use of its computations:
\emph{step-size adaptation}.
Indeed, most modern, non-probabilistic ODE solvers perform local error control with adaptive step-size selection, allowing them to compute the result up to a desired precision while avoiding unnecessary computational work.
In the following, we first review error estimation and step-size control in classic numerical solvers, then provide a principled objective for error control of probabilistic ODE solvers and describe the full step-size adaptation algorithm.

\subsection{Error Control in Classic Solvers}
An important aspect of numerical analysis is to monitor and control the error of a method.
Commonly, a distinction is made between two kinds of errors:
The \emph{local} error describes the error that the algorithm introduces after a single step,
whereas the \emph{global} error is the cumulative error of the computed solution caused by multiple iterations.
The global error is typically not of practical interest for error monitoring and control
\citep[Chapter II.4]{hairer2008solving},
and instead the local error is estimated, for example through
Richardson extrapolation
\citep[Theorem 4.1]{hairer2008solving},
or more commonly via embedded Runge-Kutta methods
\citep[Chapter II.4]{hairer2008solving}
or the Milne device
\citep{10.1145/355626.355636}.
These estimates are then used in the step-size control algorithm, which ensures that the chosen step sizes are sufficiently small to yield the desired precision of the computed result, while being sufficiently large to avoid unnecessary computational work.
Common control algorithms for step-size selection include proportional control
\citep[Chapter II.4]{hairer2008solving}
and proportional-integral (PI) control
\citep{Gustafsson1988}.

\subsection{Error Control in Probabilistic Solvers}
In Gaussian filtering, the natural object to consider for error estimation and control are the residuals \({(z_n - \hat{z}_n)}\).
\citet{schober16_probab_model_numer_solut_initial_value_probl}
show
how to use this quantity for both uncertainty calibration
(as presented in \cref{sec:uc-schober}) and local error control.
We generalize their error control objective to be applicable to all presented algorithms and uncertainty calibration methods.

After calibration, the extended Kalman filtering algorithm approximates
(see \cref{eq:update})
\begin{align}
  p(z_{n}\,|\,z_{1:n-1}) &\approx \mathcal{N}\left(z_{n}; \hat{z}_n, S_n \right),
\end{align}
with
\begin{subequations}
\begin{align}
  \hat{z}_n &= E_1 \mu_n^P - f\left( E_0 \mu_n^P, t_n \right), \\
  S_n &= H_n \left( A_n \Sigma_{n-1}^F A_n\T + \left( \breve{Q}_n \otimes \hat{\Gamma}_n \right) \right) H_n\T,
\end{align}
\end{subequations}
where \(H_n\) can correspond to either the zeroth or the first order linearization, and \(\hat{\Gamma}_n\) has been estimated through one of the approaches of \cref{sec:uncertaintycalibration}.

For step-size adaptation we want to control \emph{local} errors, and therefore assume an error-free solution estimate at time \(t_{n-1}\).
With \(\Sigma_{n-1}^F = 0\) we
obtain the approximation
\begin{equation}
  \begin{split}
    p( (z_{n} - \hat{z}_n) &\,|\, z_{1:n-1}) \approx \\
    &\mathcal{N}\left(z_{n} - \hat{z}_n; 0 , H_n \left( \breve{Q}_n \otimes \hat{\Gamma}_n \right) H_n\T \right).
  \end{split}
\end{equation}

Finally, we define the objective for local error control \(D_n \in \mathbb{R}^d\) as the (local) standard deviations of the residual vector \((z_{n}-\hat{z}_n)\), given by
\begin{equation}
  \label{eq:local-error}
  \left( D_n \right)_i := \left(H_n \left( \breve{Q}_n \otimes \hat{\Gamma}_n \right) H_n\T \right)_{ii}^{1/2}, \quad i \in \{1, \dots, d\}.
\end{equation}

For the EKF0 algorithm and time-varying, scalar \({\Gamma=\sigma^2 I_d}\), we recover the expected error used by
\citet{schober16_probab_model_numer_solut_initial_value_probl}
for step-size control (see also \citet{10.1145/355626.355636}),
but the general formulation in \cref{eq:local-error} can also be used with on-line quasi-MLE for time-fixed \(\Gamma\) (\cref{sec:uc-fixedmle}) and in combination with the EKF1.

\begin{figure*}[!t]%
  \centering%
  \begin{subfigure}[t]{0.32\linewidth}
    \centering
    \includegraphics[width=\linewidth]{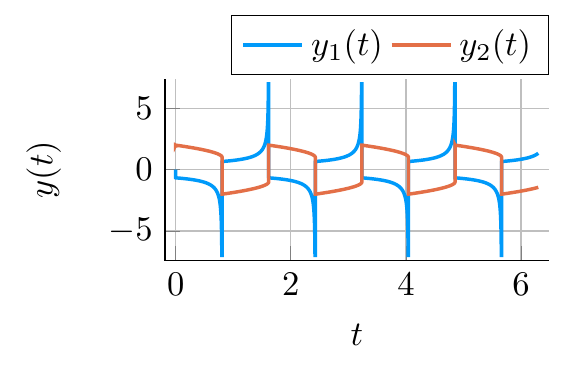}%
  \end{subfigure}
  ~
  \begin{subfigure}[t]{0.32\linewidth}
    \centering
    \includegraphics[width=\linewidth]{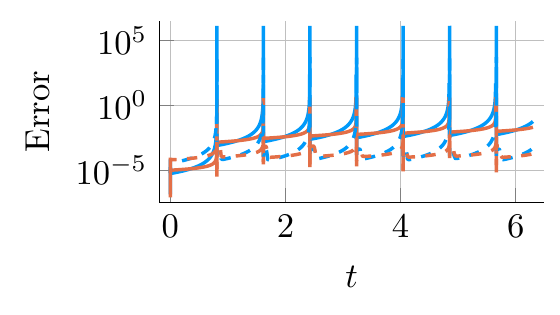}%
  \end{subfigure}
  ~
  \begin{subfigure}[t]{0.32\linewidth}
    \centering
    \includegraphics[width=\linewidth]{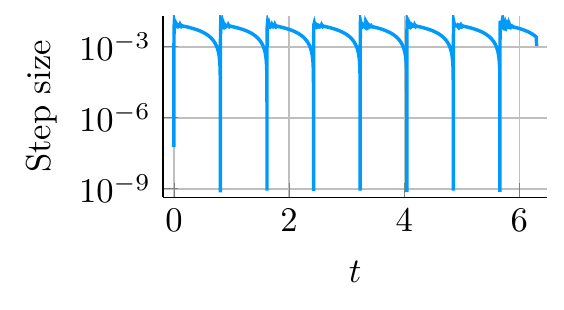}%
  \end{subfigure}
  \caption{
    \emph{Probabilistic solution with step-size adaptation of the Van der Pol equations}.
    Left:
    Mean of the posterior distribution over solutions, showing only values in the $(-7,7)$ interval.
    Middle:
    Absolute errors (solid lines), and standard deviations of the posterior marginals as error estimates (dashed lines), shown in log-scale.
    Right:
    Step sizes of accepted steps throughout the solve.
    During the stiff phases, the step sizes get drastically decreased by the step-size controller.
  }
  \label{fig:vdp}
\end{figure*}

\subsection{Step-Size Selection}
Following \citet[Chapter II.4]{hairer2008solving},
the step-size controller
aims to select step sizes, as large as possible, but while satisfying componentwise,
for \({i \in \{1, \dots, d\}}\),
\begin{equation}
  \label{eq:stepsize-objective}
  (D_n)_i \leq \varepsilon_i, \quad
  \varepsilon_i := \tau_\text{abs} + \tau_\text{rel} \cdot \max \left( \left| (\hat{y}_{n-1})_i \right|, \left| (\hat{y}_{n})_i \right|  \right),
\end{equation}
where \(\tau_\text{abs}\) and \(\tau_\text{rel}\) are the prescribed absolute and relative tolerances, respectively, and
\(\hat{y}_{n-1} := E_0 \mu_{n-1}^F\), \(\hat{y}_{n} := E_0 \mu_n^F\) are solution estimates of the numerical solver.
To do so, we define the following measure of error as control objective,
\begin{equation}
  E := \sqrt{\frac{1}{d} \sum_{i=1}^d \left( \frac{ \left( D_n \right)_i }{\varepsilon_i} \right)^2}.
\end{equation}

The proportional control algorithm
compares the control objective \(E\) to \(1\) to find the optimal step size.
If \(E \leq 1\) holds, the computed step is accepted and the integration continues.
Otherwise, the step is rejected as too inaccurate and is repeated.
In both cases, a new step size which will likely satisfy \cref{eq:stepsize-objective} is computed as
\begin{equation}
  h_\text{new} = h \cdot \rho \left( \frac{1}{E} \right)^{\frac{1}{q+1}},
\end{equation}
making use of the local convergence rate of \(q+1\),
as used by
\citet{schober16_probab_model_numer_solut_initial_value_probl}
and shown by
\citet{Kersting2020}.
The parameter \(\rho \in (0, 1]\) is a safety factor to increase the probability that the next step will be acceptable, and we additionally limit the rate of change
\(\eta_\text{min} \leq h_{n+1}/h_n \leq \eta_\text{max}\)
\citep{hairer2008solving}.
In our experiments, we set
\(\rho := 0.9\), \(\eta_{min} := 0.2\), \(\eta_{max} := 10\).

The formulation of the local error control objective in
\cref{eq:local-error}
also lends itself to other control algorithms.
Notably, PI control
\citep{Gustafsson1988}
can be an interesting alternative to proportional control when applied to mildly stiff problems,
and has been successfully applied for the related class of Nordsieck methods by
\citet{Bras_2013}.

\section{RELATED WORK}
Our contribution fits in the formulation of IVPs as problems of Bayesian state estimation
\citep{DBLP:journals/sac/TronarpKSH19,tronarp20_bayes_ode_solver}.
By using Gaussian filtering methods, these solvers are able to efficiently compute posterior distributions over solutions
\citep{DBLP:conf/uai/KerstingH16},
which converge to the true solution at high polynomial rates
\citep{Kersting2020}.
In practice, the calibration of the posterior uncertainties depends on specific model hyperparameters.
This paper reviews and extends previously proposed global and local calibration methods
\citep{DBLP:journals/sac/TronarpKSH19,schober16_probab_model_numer_solut_initial_value_probl}.
The presented step-size controller builds on the algorithm suggested by
\citet{schober16_probab_model_numer_solut_initial_value_probl}.

A different line of work on probabilistic numerical solvers for ODEs
aims to represent the distribution over solution with a set of sample paths
\citep{DBLP:journals/sac/ConradGSSZ17, DBLP:journals/sac/AbdulleG20, DBLP:journals/sac/LieSS19, DBLP:conf/nips/TeymurLSC18, DBLP:conf/nips/TeymurZC16, chkrebtii2016, DBLP:journals/sac/TronarpKSH19}.
While these methods are able to capture arbitrary, non-Gaussian distributions, they come at an increased computational cost.

\section{EXPERIMENTS}
\label{sec:experiments}
To evaluate the presented methodology, we provide three sets of experiments.
First, we highlight the practical necessity of step-size adaptation by solving a stiff version of the Van der Pol model.
Next, we compare the different uncertainty calibration methods presented in
\cref{sec:uncertaintycalibration}.
Finally, we assess the practical performance of the probabilistic ODE solvers by comparing to a classic Runge-Kutta 4/5 method.
The code for the implementation and experiments is publicly available on
gihub\footnote{\url{https://github.com/nathanaelbosch/capos}}.

\begin{figure*}[tb]
  \centering
  \begin{subfigure}[b]{\linewidth}
    \centering
    \includegraphics[width=\linewidth]{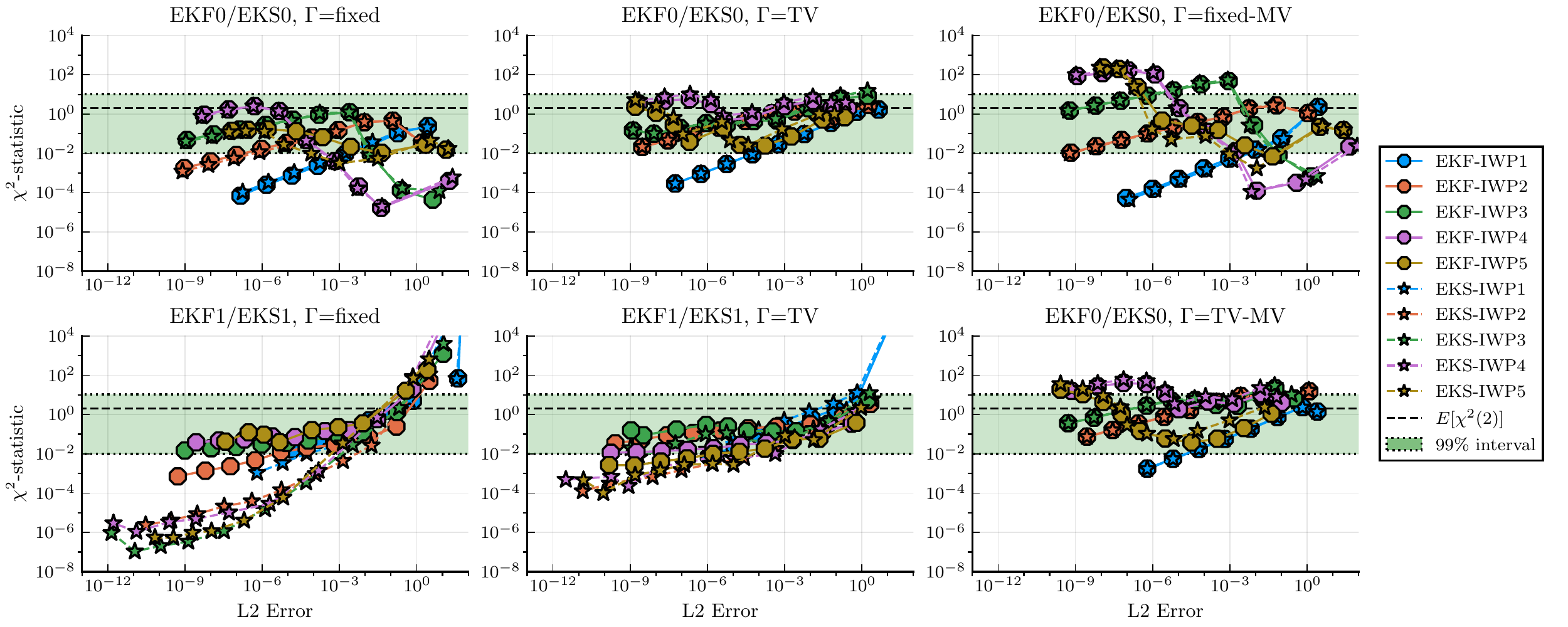}
  \end{subfigure}
  \caption{
    \emph{Uncertainty calibration across configurations}.
    In each subfigure, a specific combination of filtering algorithm (EKF0/EKS0 or EKS1/EKF1) and calibration method is evaluated, the latter including fixed and time-varying (TV) diffusion models, as well as their multivariate versions (fixed-MV, TV-MV).
    A well-calibrated solver should provide
    \(\chi^2\)-statistics inside the \(99\%\) credible interval (green).
  }
  \label{fig:calibration-comparison}
\end{figure*}

\begin{figure*}[t]%
  \centering%
  \begin{subfigure}[t]{0.35\linewidth}
    \centering
    \includegraphics[width=\linewidth]{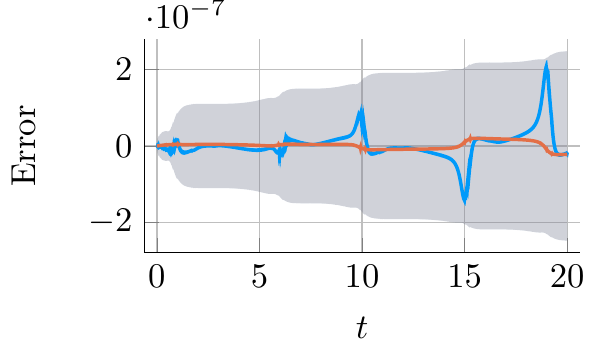}%
  \end{subfigure}
  \hspace{-2.0cm}
  \begin{subfigure}[t]{0.35\linewidth}
    \centering
    \includegraphics[width=1.525\linewidth]{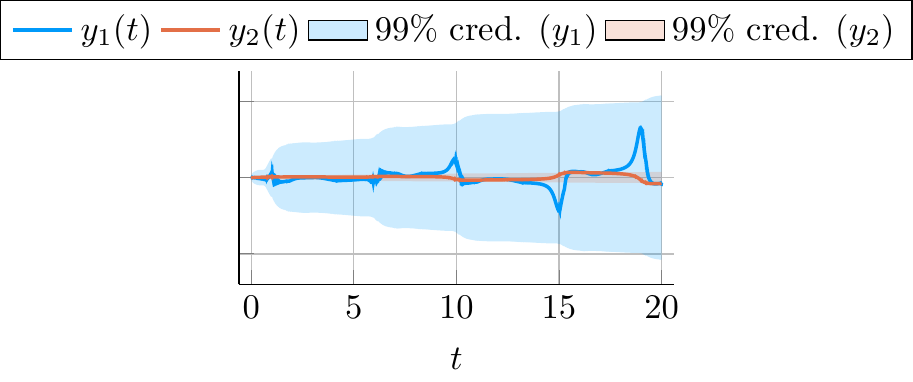}%
  \end{subfigure}
  \hspace{0.6cm}
  \begin{subfigure}[t]{0.35\linewidth}
    \centering
    \includegraphics[width=0.79\linewidth]{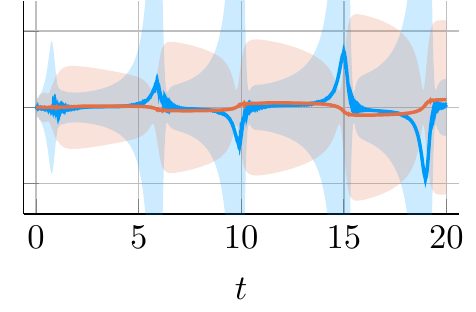}%
  \end{subfigure}
  \caption{
    \emph{Qualitative comparison of the different time-varying uncertainty models}.
    The EKS0 with scalar diffusion (left) estimates a single credible band shared accross for both dimensions, but
    the multivariate model (middle) is able to attribute different uncertainties to each dimension.
    By including information on the derivatives of the ODE, the EKS1 with scalar diffusion (right) returns structured and dynamic uncertainty estimate.
  }
  \label{fig:uncertainty-structure}
\end{figure*}

\subsection{Stiff Van der Pol}
The Van der Pol model \citep{vanderpol}
describes a non-conservative oscillator with non-linear damping, and can be written in the two-dimensional form:
\begin{equation}
  \begin{split}
  \dot{y_1} &= y_2, \\
  \dot{y_2} &= \mu \left((1-y_1^2) y_2 - y_1\right),
\end{split}
\end{equation}
with a positive stiffness constant \(\mu>0\).

To highlight the importance of step-size adaptation, and to demonstrate the A-stability of the EKS1
\citep{tronarp20_bayes_ode_solver},
we consider a very stiff version of the Van der Pol model and set \(\mu=10^6\).
We solve the IVP on the time interval \([0, 6.3]\) with initial value \(y(t_0) = [0, \sqrt{3}]\T\)
using the EKS1 algorithm
with an IWP3 prior, for
absolute and relative tolerance specified as \(10^{-6}\) and \(10^{-3}\), respectively.
We were not able to solve this same IVP with the EKS0 (which does not possess this stability property).
The reference solution has been computed with the
A-stable \(5\)th order implicit Runge-Kutta method Radau IIA
\citep{hairer1987solving},
implemented as \texttt{RadauIIA5} in the Julia DifferentialEquations.jl suite
\citep{rackauckas2017differentialequations},
for absolute and relative tolerances set to \(10^{-14}\).

\Cref{fig:vdp} shows the results, including the solution posterior, errors and error estimates, and the step sizes during the solve.
The posterior mean returned by the solver achieved a final error of
\(\|\hat{y}(T) - y^*(T)\| = 6.17 \times 10^{-2}\), where \(y^*\) denotes the
reference solution.
The step sizes shown in \cref{fig:vdp} (right) change drastically during the
solve, and decrease from values \(h \sim 10^{-2}\) down to \( h < 10^{-8} \) during the stiff phases.
If one were to solve this IVP without step-size control, that is, with EKS1 and an IWP3 prior but with fixed steps of size \(10^{-8}\), one would perform \(6.3 \times 10^{8}\) solver iterations.
Compared to the \(\sim 1\) second runtime the adaptive step method required for its \(23824\) iterations (out of which \(6977\) steps were rejected), the fixed-step solver would require more than \(7\) hours.
This demonstrates the importance of adaptive step-size control for efficient use of computational resources.

\subsection{Comparison of Calibration Methods}
\label{sec:experiments:calibration}

We evaluate the various algorithms and calibration methods on the FitzHugh-Nagumo model, given by the ODE
\begin{equation}
  \begin{split}
    \dot{y}_1 &= c \left( y_1 - \frac{y_1^3}{3} + y_2 \right), \\
    \dot{y}_2 &= -\frac{1}{c} \left( y_1 -  a - b y_2 \right),
  \end{split}
\end{equation}
with parameters
(\(a=0.2\), \(b=0.2\), \(c=3.0\)),
initial value
\(y(0) = [-1, 1]\T\),
and time span \([0, 20]\).

In the first part of this experiment, the uncertainty calibration is evaluated across a large range of solver configurations and tolerances
(\(\tau_\text{abs} = 10^{-4}, \dots, 10^{-13}\),
\(\tau_\text{rel} = 10^{-1}, \dots, 10^{-10}\)).
To assess the quality of the uncertainty calibration, we use the \(\chi^2\)-statistics
\citep{bar2004estimation}
defined by
\begin{align}
  \chi^2 = \frac{1}{N} \sum_{i=1}^N r(t_i)\T \Cov(y(t_i))^{-1} r(t_i),
\end{align}
where \(r(t_i) := (y^*(t_i) - \mathbb{E}[y(t_i)])\) are the residuals and
\(\mathbb{E}[y(t_i)]\) and \(\Cov(y(t_i))\) are computed on the posterior distribution returned by the probabilistic solver.
A well calibrated model achieves \(\chi^2 \approx d\).
If \(\chi^2 < d\) or
\(\chi^2 > d\) we refer to the solution as
\emph{underconfident}
or
\emph{overconfident},
respectively.

\Cref{fig:calibration-comparison}
visualizes the full comparison of the uncertainty calibration methods and suggests some empirical findings.
For most configurations the time-varying uncertainty models seem to be better calibrated than the fixed models.
For a fixed choice of algorithm (e.g. EKF0) and order (e.g. IWP5) their calibration varies less across the different tolerance levels.
The multi-variate models seem to not have a large impact on the \(\chi^2\)-statistics, both for the fixed and time-varying models.
The first-order linearization approaches EKF1/EKS1 tend to become underconfident, but achieve the lowest errors.

To complement this summarized evaluation based on the \(\chi^2\) statistics, we visualize the qualitative behaviour of the different time-varying uncertainty models in
\cref{fig:uncertainty-structure}.
To be comparable, all models share an IWP3 prior and tolerances \(\tau_\text{abs}=10^{-10}\), \(\tau_\text{rel}=10^{-7}\).
The scalar, time-varying (TV) approach attributes the same credible bands to both
dimensions (shown in grey), whereas the
multivariate, time-varying (TV-MV) model is able to lift this restriction and
estimates a large uncertainty for the first dimension (blue) and barely visible uncertainties for the second dimension (orange).
However, only the first order linearization of the EKS1 seems to properly describe the structural properties of the true solution in its posterior estimate.

For experiments on additional problems, including classic work-precision diagrams for all methods, see \cref{sec:supplementary:calibration}

\subsection{Comparison with Dormand--Prince 4/5}
\label{sec:dp5-comparison}
This experiment assesses the performance of the developed methodology and compares the probabilistic solvers of \(5\)th order to the classic, widely used Runge-Kutta 4/5 method by
\citet{dormand1980family},
implemented as \texttt{DP5} in the Julia DifferentialEquations.jl suite
\citep{rackauckas2017differentialequations}.
The comparison was made on the Lotka-Volterra equations
\citep{lotka,volterra},
which describe the dynamics of biological systems in which two species interact, one as a predator and the other as prey.
The IVP is given by the ODE
\begin{equation}
  \begin{split}
    \dot{y_1} &= \alpha y_1 - \beta y_1 y_2, \\
    \dot{y_2} &= -\gamma y_1 + \delta y_1 y_2,
  \end{split}
\end{equation}
with initial value \(y(0) = [1,1]\T\) and parameters (\(\alpha=1.5\), \(\beta=1\), \(\gamma=3\), \(\delta=1\)), on the time span \([0,10]\).

\begin{figure}[tb]
  \centering
  \includegraphics[width=\linewidth]{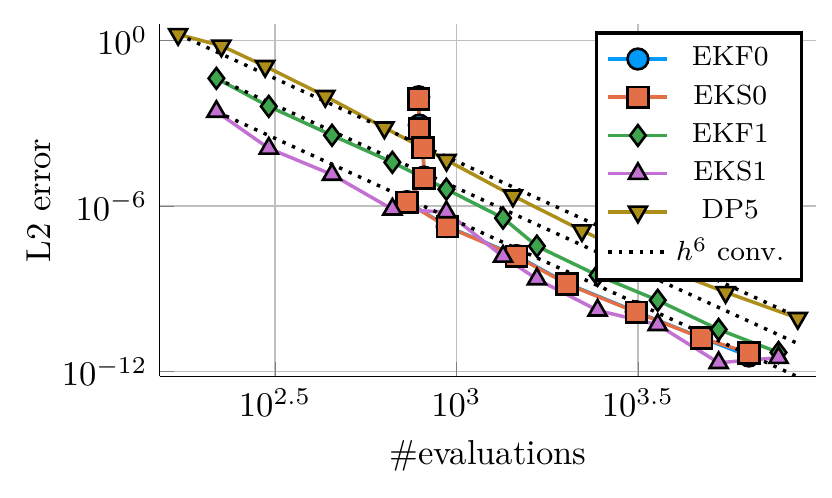}%
  \caption{
    \emph{Comparison to Dormand--Prince 4/5 (DP5)}.
    All probabilistic ODE solvers use an IWP-5 prior and a scalar, time-varying diffusion model.
    The comparison is made on the Lotka-Volterra equations.
    The EKF0 and EKS0 performed similarly and are therefore difficult to separate visually in the figure.
  }
  \label{fig:dp5-comparison-lv}
\end{figure}

\Cref{fig:dp5-comparison-lv} shows the results in a work-precision diagram.
We observe convergence rates of order \(6\) for all methods, one order higher than the expected global convergence rate of order \(5\) for the DP5 algorithm
\citep{hairer2008solving}
and for Gaussian ODE filters with IWP-\(5\) prior
\citep{Kersting2020,tronarp20_bayes_ode_solver}.
It can also be seen that the EKF0 requires more function evaluations than expected for high-tolerance settings.
Out of all compared methods, the EKS1 seems to require the least number of evaluations of the function and its Jacobian to achieve a specified error, matching the performance of the EKF0 and EKS0 while demonstrating a stable behaviour for high tolerances.

Similar work-precision diagrams for these methods on additional problems are provided in \cref{sec:supplementary:dp5}.

\section{CONCLUSION}
In this paper, we introduced and discussed various models and methods for uncertainty calibration in Gaussian ODE filters, and presented parameter estimates for both fixed and time-varying, as well as scalar and multivariate diffusion models.
The probabilistic observation model of these methods provides a calibrated objective for local error control, enabling the implementation of classic step-size selection algorithms.

The resulting, efficiently computable posteriors have been empirically evaluated for a wide range of tolerance levels, demonstrating decent error calibration in particular for the time-varying diffusion models.
Of all compared methods, the first-order linearization of the EKS1 seems to provide the most expressive posterior covariances, while also efficiently computing accurate solutions -- requiring, in our benchmarks, less evaluations than the well-known Dormand--Prince 4/5 method to reach a specified tolerance level.

\subsubsection*{Acknowledgements}
The authors gratefully acknowledge financial support by the German Federal Ministry of Education and Research (BMBF) through Project ADIMEM (FKZ 01IS18052B), and
financial support by the European Research Council through ERC StG Action 757275 / PANAMA; the DFG Cluster of Excellence “Machine Learning - New Perspectives for Science”, EXC 2064/1, project number 390727645; the German Federal Ministry of Education and Research (BMBF) through the Tübingen AI Center (FKZ: 01IS18039A); and funds from the Ministry of Science, Research and Arts of the State of Baden-Württemberg.
The authors also thank the International Max Planck Research School for Intelligent Systems (IMPRS-IS) for supporting N. Bosch.

The authors are grateful to Nicholas Krämer for many valuable discussions.
They further thank Hans Kersting, Jonathan Wenger and Agustinus Kristiadi for helpful feedback on the manuscript.

\bibliography{references}

\clearpage

\appendix
\onecolumn
\aistatstitle{Calibrated Adaptive Probabilistic ODE Solvers \\
  Supplementary Materials}

\section{Proof of \cref{prop:fixedMV}}
\label{sec:proofs}
\begin{proof}
\label{sec:proofs:ekf0-fixed}
The proof is structured as follows.
First, we show by induction that an initial covariance
\(\Sigma_0 = \breve{\Sigma}_0 \otimes \Gamma\)
implies covariances
\(\Sigma_{n}^P = \breve{\Sigma}_{n}^p \otimes \Gamma \),
\(\Sigma_{n}^F = \breve{\Sigma}_{n}^F \otimes \Gamma\),
and
\(S_n = \breve{S}_n \cdot \Gamma\),
for all \(n\).
Then, for measurement covariances \(S_n\) of such form, we can compute the (quasi) maximum likelihood estimate
\(\hat{\Gamma}\).

Assume \(\Sigma_{n-1}^F = \breve{\Sigma}_{n-1}^F \otimes \Gamma\).
Using the mixed product property and the associativity of the Kronecker product, the covariance of the prediction \(\Sigma_n^P\) can be written as
\begin{align*}
  \Sigma_n^P
  = A_n \Sigma_{n-1}^F A_n\T + Q_n
  = \left(\breve{A}_n \otimes I_d \right) \left(\breve{\Sigma}_{n-1}^F \otimes \Gamma \right) \left(\breve{A}_n \otimes I_d \right) \T + \left(\breve{Q}_n \otimes \Gamma \right)
  = \left(\breve{A}_n \breve{\Sigma}_{n-1}^F \breve{A}_n\T + \breve{Q}_n \right) \otimes \Gamma
  = \breve{\Sigma}_n^P \otimes \Gamma,
\end{align*}
where
\(\breve{\Sigma}_n^P := \breve{A}_n \breve{\Sigma}_{n-1}^F \breve{A}_n\T + \breve{Q}_n\).
Next, using \(H_n = e_1\T \otimes I_d\) (EKF0), the measurement covariance \(S_n\) is given by
\begin{align*}
  S_n
  = H_n \Sigma_n^P H_n\T
  = \left(e_1\T \otimes I_d \right) \left(\breve{\Sigma}_n^P \otimes \Gamma \right) \left(e_1\T \otimes I_d \right)\T
  = \left(e_1\T \breve{\Sigma}_n^P e_1 \right) \otimes \Gamma
  = \breve{S}_n \cdot \Gamma.
\end{align*}
Finally,
the filtering covariance can be computed with the update step \cref{eq:update} as
\begin{align*}
  \Sigma_n^F
  &= \Sigma_n^P - K_n S_n K_n \\
  &= \Sigma_n^P - \Sigma_n^P H_n\T S_n^{-1} H_n \Sigma_n^P \\
  &= \left(\breve{\Sigma}_n^P \otimes \Gamma \right) - \left(\breve{\Sigma}_n^P \otimes \Gamma \right) \left(e_1\T \otimes I_d \right)\T \left(\breve{S}_n \otimes \Gamma \right)^{-1} \left(e_1\T \otimes I_d \right) \left(\breve{\Sigma}_n^P \otimes \Gamma \right) \\
  &= \left(\breve{\Sigma}_n^P  - \breve{\Sigma}_n^P e_1 \breve{S}_n^{-1} e_1\T \breve{\Sigma}_n^P \right) \otimes \Gamma \\
  &= \breve{\Sigma}_n^F \otimes \Gamma,
\end{align*}
with
\(\breve{\Sigma}_n^P := \breve{\Sigma}_n^P  - \breve{\Sigma}_n^P e_1 \breve{S}_n^{-1} e_1\T \breve{\Sigma}_n^P \).
This concludes the first part of the proof.

It is left to compute the (quasi) MLE by maximizing the log-likelihood
\({\log p(z_{1:N}) = \log \prod_{n=1}^N \mathcal{N} \left( z_n; \hat{z}_n, S_n \right)}\).
For diagonal \(\Gamma\) and zero measurements \(z_n = 0\), we obtain
\begin{align*}
  \hat{\Gamma}
  &= \argmax_\Gamma \log p(z_{1:N}) \\
  &= \argmax_\Gamma \sum_{n=1}^N \log \mathcal{N} \left( 0; \hat{z}_n, \breve{S}_n \cdot \Gamma \right) \\
  &= \argmax_\Gamma \sum_{n=1}^N \left( - \frac{1}{2} \log |\breve{S}_n \cdot \Gamma| - \frac{1}{2} \hat{z}_n\T \left( \breve{S}_n \cdot \Gamma \right)^{-1} \hat{z}_n \right) \\
  &= \argmax_\Gamma \sum_{n=1}^N \left( - \frac{\log \breve{S}_n^d + \sum_{i=1}^d \log \Gamma_{ii}}{2} - \sum_{i=1}^d \frac{(\hat{z}_n)_i^2}{2 \breve{S}_n \Gamma_{ii}} \right) \\
  &= \argmax_\Gamma \sum_{i=1}^d \left( - \frac{N \log \Gamma_{ii} }{2} - \sum_{n=1}^N \frac{(\hat{z}_n)_d^2}{2 \breve{S}_n \Gamma_{ii}} \right).
\end{align*}
Each diagonal element
\(\hat{\Gamma}_{ii}\) can be computed independently by taking the derivative setting it to zero:
\begin{align*}
  0 = - \frac{N}{2 \Gamma_{ii}} + \sum_{n=1}^N \frac{(\hat{z}_n)_d^2}{2 \breve{S}_n \Gamma_{ii}^2} , \qquad \forall i \in \{1, \dots, d\}.
\end{align*}
The solution to this equation provides the (quasi) maximum likelihood estimate for \(\hat{\Gamma}\)
\begin{align*}
  \hat{\Gamma}_{ii} = \frac{1}{N} \sum_{n=1}^N \frac{(\hat{z}_n)_i^2}{\breve{S}_n} , \qquad \forall i \in \{1, \dots, d\}.
\end{align*}
\end{proof}

\section{Additional Experiments}
These experiments include two additional IVPs.

\paragraph{Logistic Equation}
We consider an initial value problem, given by the logistic equation
\begin{equation}
  \dot{y}(t) = r y(t) (1-y(t)),
\end{equation}
with parameter \(r=3\),
integration interval \([0, 2.5]\),
and initial value \(y(0) = 0.1\).
Its exact solution is given by
\begin{equation}
  y^*(t) = \frac{\exp(rt)}{1/y_0 - 1 + \exp(rt)}.
\end{equation}

\paragraph{Brusselator}
The Brusselator is a model for multi-molecular chemical reactions, given by the ODEs
\begin{equation}
  \begin{split}
    \dot{y_1} &= 1 + y_1^2 y_2 - 4 y_1, \\
    \dot{y_2} &= 3 y_1 - y_1^2 y_2.
  \end{split}
\end{equation}
We consider an IVP with initial value \(y(0) = [1.5, 3]\)
on the time span \([0, 10]\).

\subsection{Performance and Calibration on Additional Problems}
\label{sec:supplementary:calibration}

\begin{figure}[tb]
  \centering
  \begin{subfigure}[b]{\linewidth}
    \centering
    \includegraphics[width=1\linewidth]{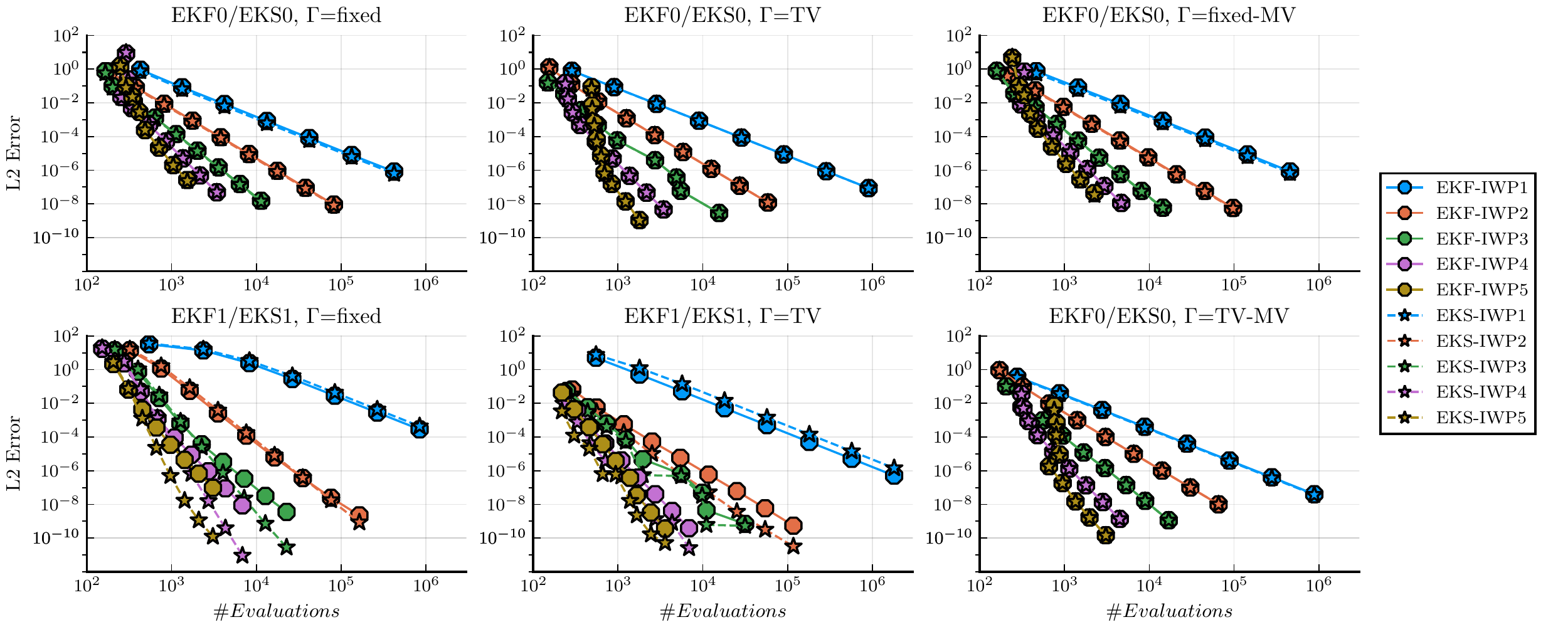}
    \caption{Work-precision diagrams.}
  \end{subfigure}

  \vspace{1em}

  \begin{subfigure}[b]{\linewidth}
    \centering
    \includegraphics[width=1\linewidth]{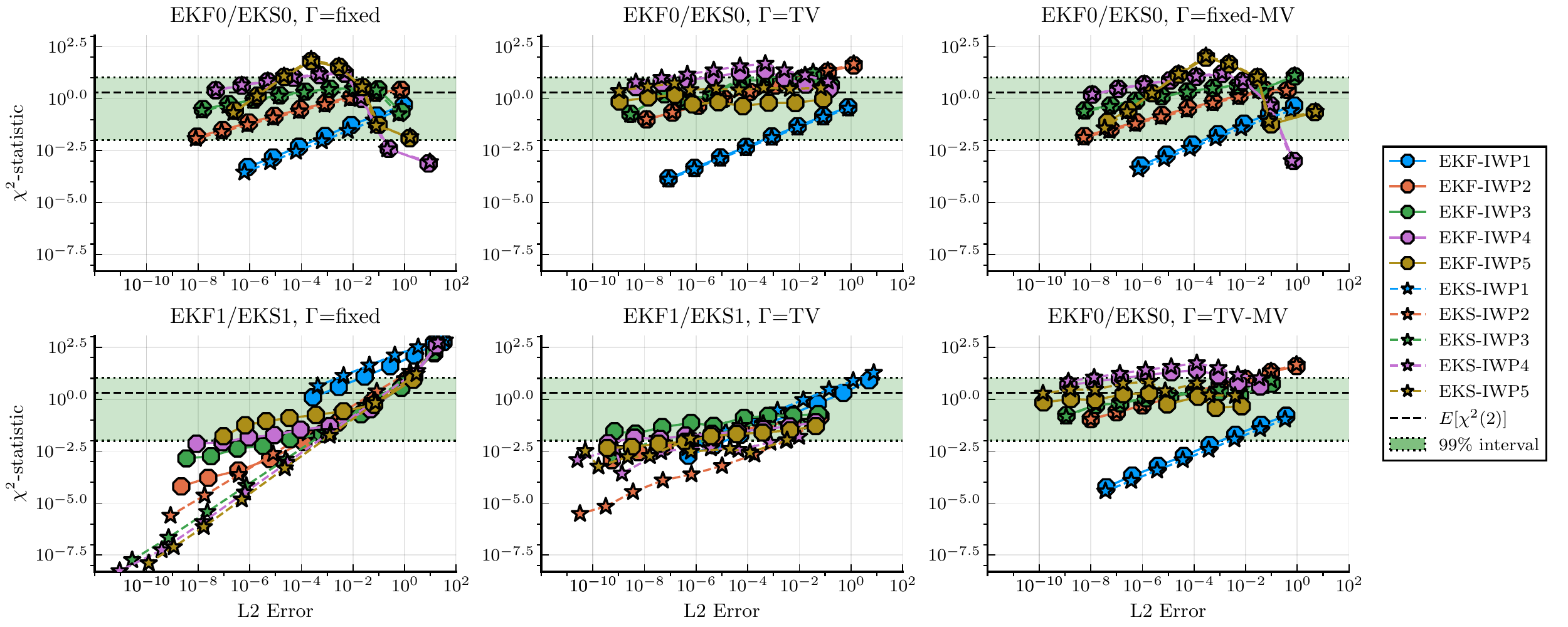}
    \caption{Uncertainty Calibration.}
  \end{subfigure}
  \caption{
    \emph{Accuracy and uncertainty calibration across configurations on the Lotka-Volterra equations}.
    In each subfigure, a specific combination of filtering algorithm (EKF0/EKS0 or EKS1/EKF1) and calibration method is evaluated, the latter including fixed and time-varying (TV) diffusion models, as well as their multivariate versions (fixed-MV, TV-MV).
  }
  \label{fig:calibration:lotkavolterra}
\end{figure}

\begin{figure}[tb]
  \centering
  \begin{subfigure}[b]{\linewidth}
    \centering
    \includegraphics[width=\linewidth]{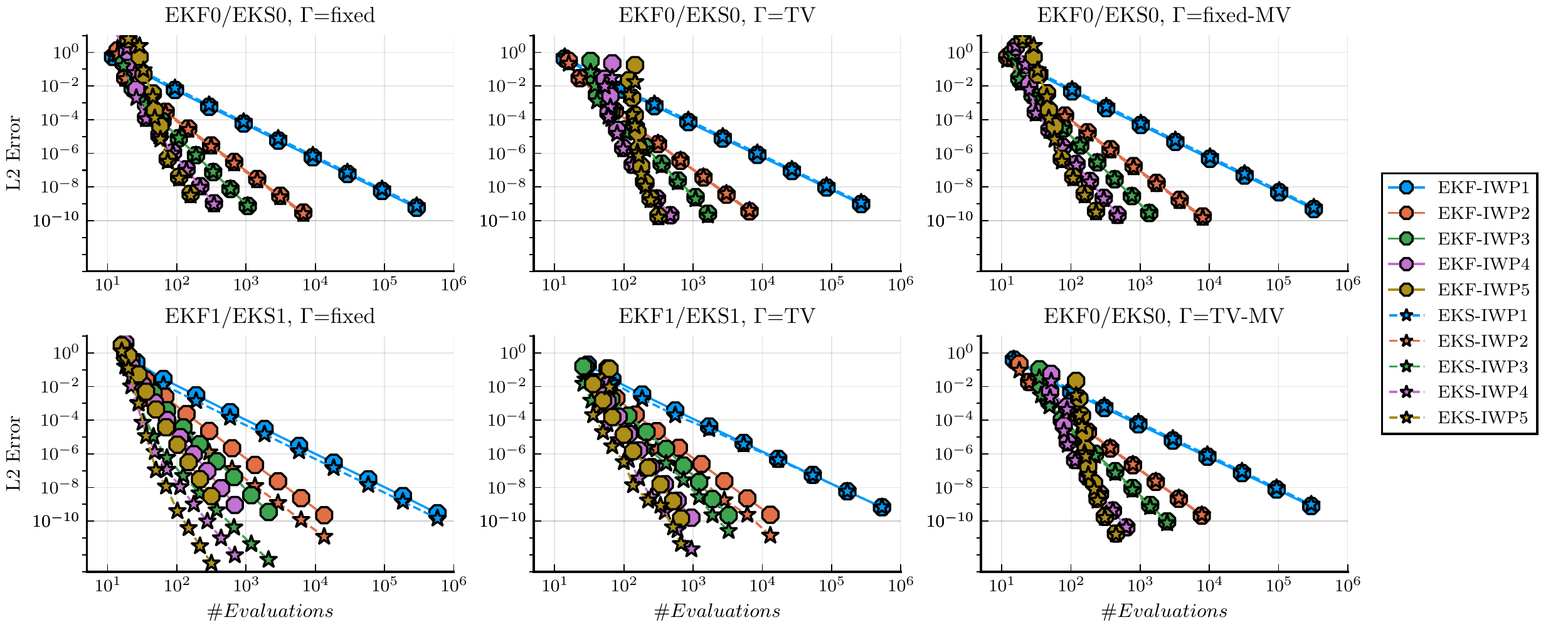}
    \caption{Work-precision diagrams.}
  \end{subfigure}

  \vspace{1em}

  \begin{subfigure}[b]{\linewidth}
    \centering
    \includegraphics[width=1\linewidth]{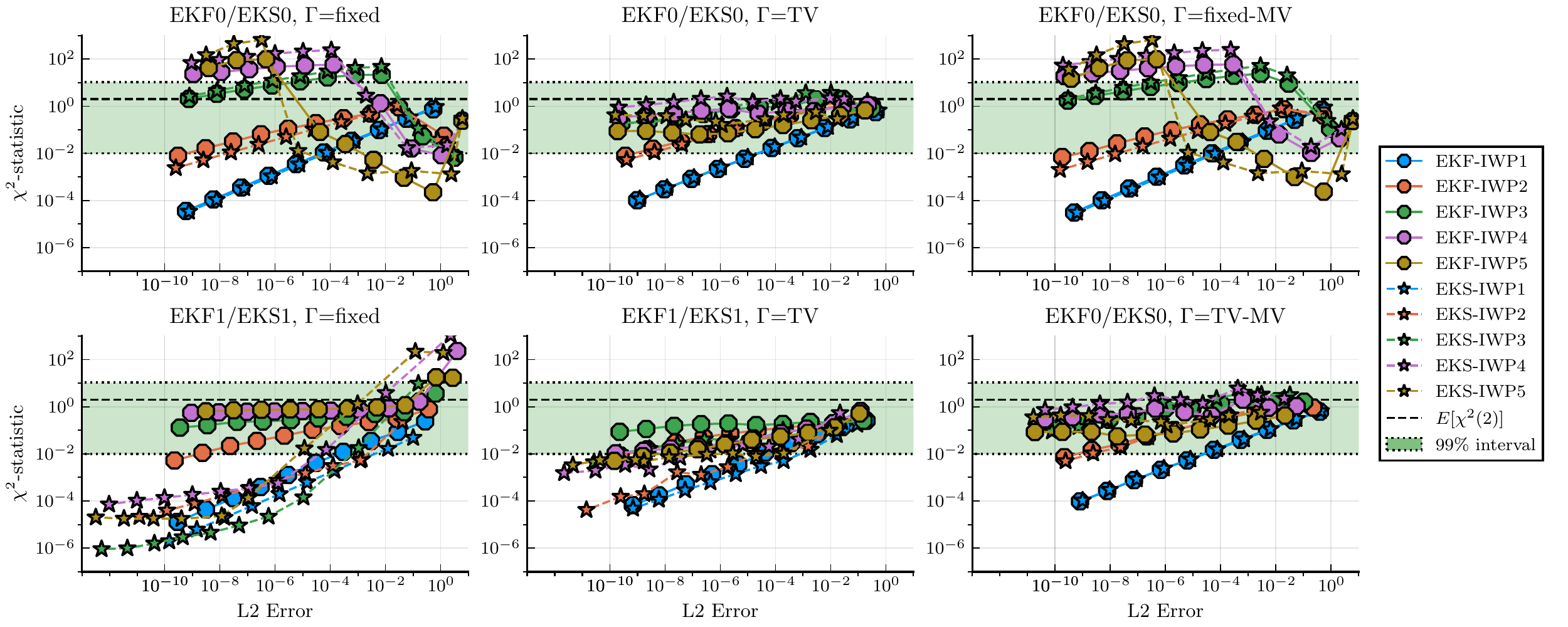}
    \caption{Uncertainty Calibration.}
  \end{subfigure}
  \caption{
    \emph{Accuracy and uncertainty calibration across configurations on the logistic equation}.
    In each subfigure, a specific combination of filtering algorithm (EKF0/EKS0 or EKS1/EKF1) and calibration method is evaluated, the latter including fixed and time-varying (TV) diffusion models, as well as their multivariate versions (fixed-MV, TV-MV).
  }
  \label{fig:calibration:logistic}
\end{figure}

\begin{figure}[tb]
  \centering
  \begin{subfigure}[b]{\linewidth}
    \centering
    \includegraphics[width=1\linewidth]{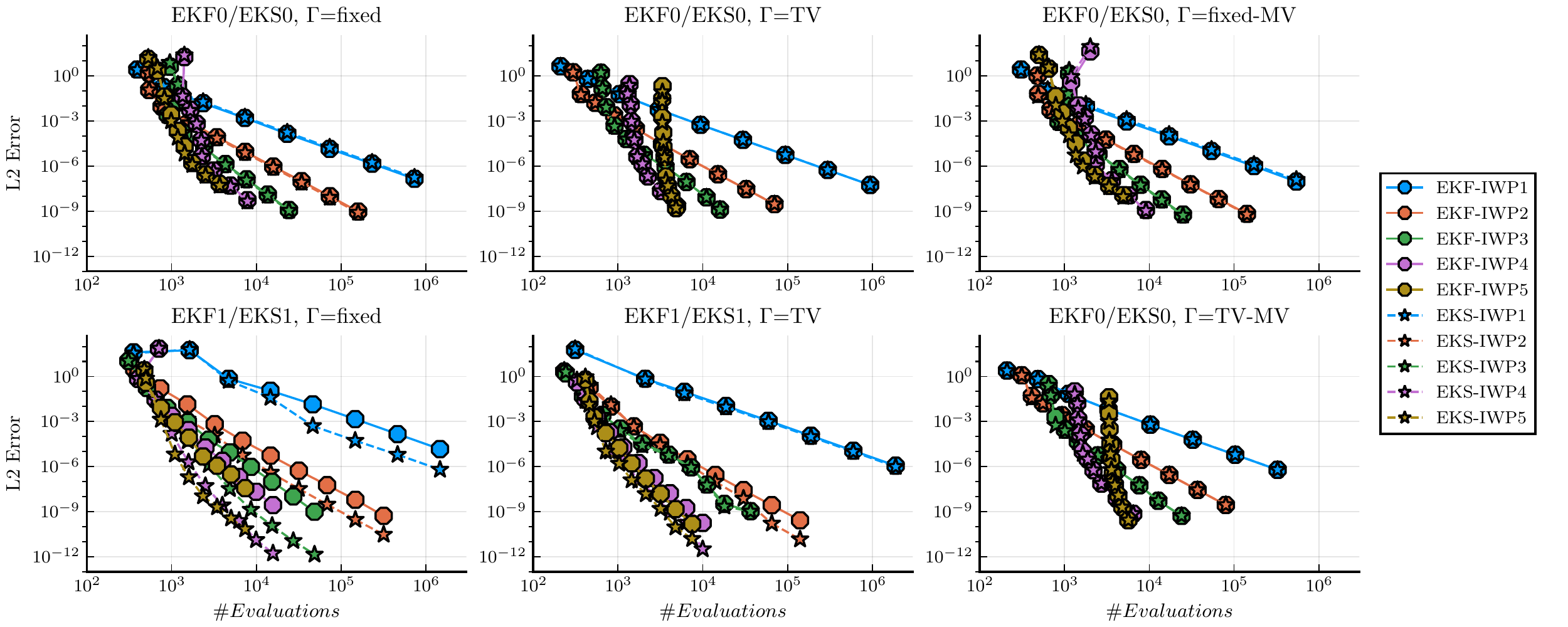}
    \caption{Work-precision diagrams.}
  \end{subfigure}

  \vspace{1em}

  \begin{subfigure}[b]{\linewidth}
    \centering
    \includegraphics[width=1\linewidth]{figures/fitzhughnagumo_calibration.pdf}
    \caption{Uncertainty Calibration.}
  \end{subfigure}
  \caption{
    \emph{Accuracy and uncertainty calibration across configurations on the FitzHugh-Nagumo equations}.
    In each subfigure, a specific combination of filtering algorithm (EKF0/EKS0 or EKS1/EKF1) and calibration method is evaluated, the latter including fixed and time-varying (TV) diffusion models, as well as their multivariate versions (fixed-MV, TV-MV).
  }
  \label{fig:calibration:fitzhughnagumo}
\end{figure}

\begin{figure}[tb]
  \centering
  \begin{subfigure}[b]{\linewidth}
    \centering
    \includegraphics[width=1\linewidth]{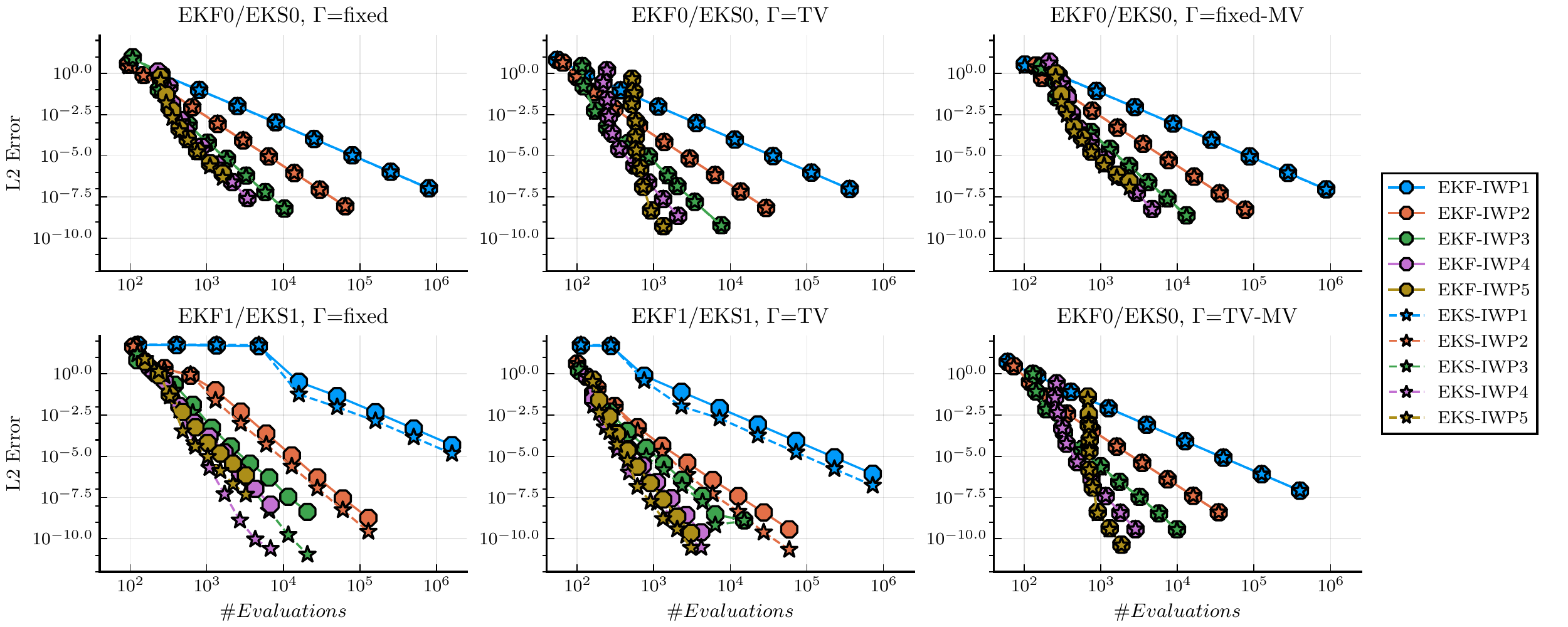}
    \caption{Work-precision diagrams.}
  \end{subfigure}

  \vspace{1em}

  \begin{subfigure}[b]{\linewidth}
    \centering
    \includegraphics[width=1\linewidth]{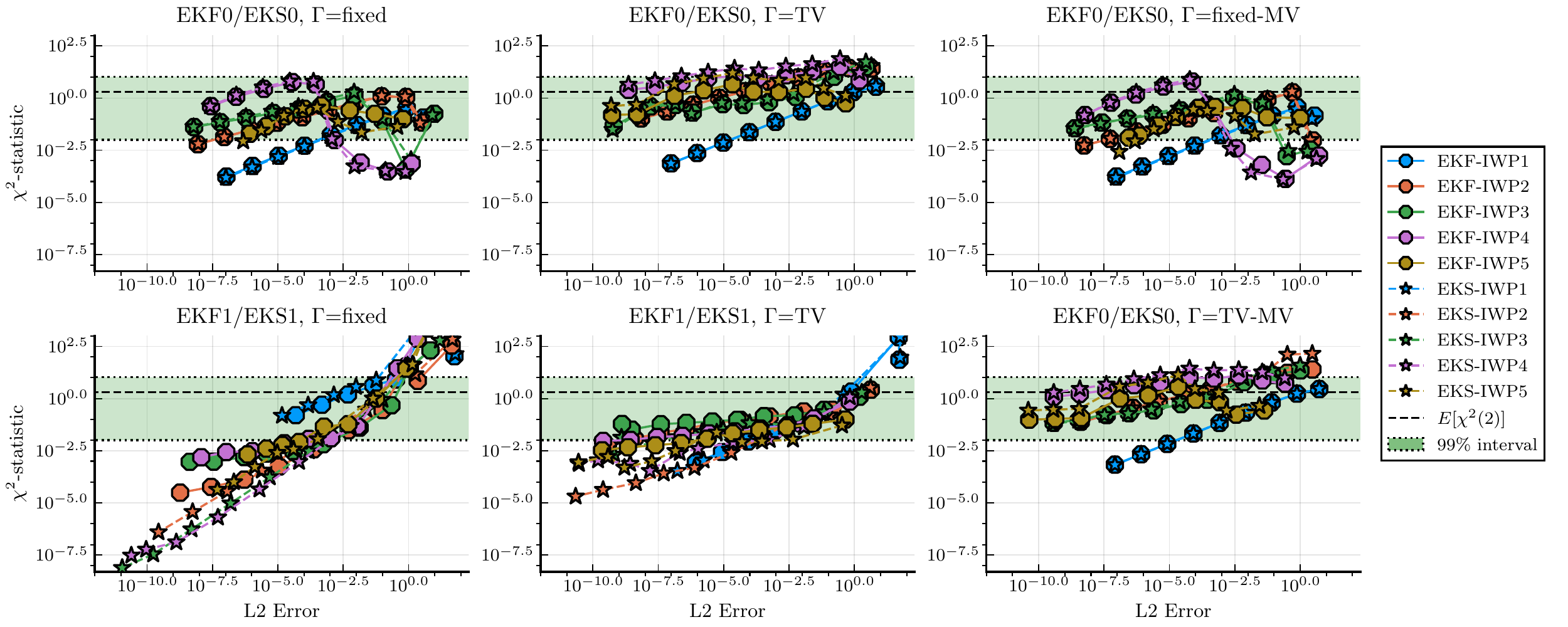}
    \caption{Uncertainty Calibration.}
  \end{subfigure}
  \caption{
    \emph{Accuracy and uncertainty calibration across configurations on the Brusselator equations}.
    In each subfigure, a specific combination of filtering algorithm (EKF0/EKS0 or EKS1/EKF1) and calibration method is evaluated, the latter including fixed and time-varying (TV) diffusion models, as well as their multivariate versions (fixed-MV, TV-MV).
  }
  \label{fig:calibration:brusselator}
\end{figure}

This section extends the results of
\cref{sec:experiments:calibration}
with evaluations on additional problems:
Lotka-Volterra (\cref{fig:calibration:lotkavolterra}),
logistic equation (\cref{fig:calibration:logistic}),
FitzHugh-Nagumo (\cref{fig:calibration:fitzhughnagumo}),
and Brusselator (\cref{fig:calibration:brusselator}).
In addition, all figures contain classic work-precision diagrams which visualize the relation of achieved error and number of evaluations, where both the evaluations of the vector field and of its Jacobian are counted.

\subsection{Comparison to Dormand-Prince 4/5}
\begin{figure}[tb]
  \centering
  \begin{subfigure}{0.48\linewidth}
    \includegraphics[width=\linewidth]{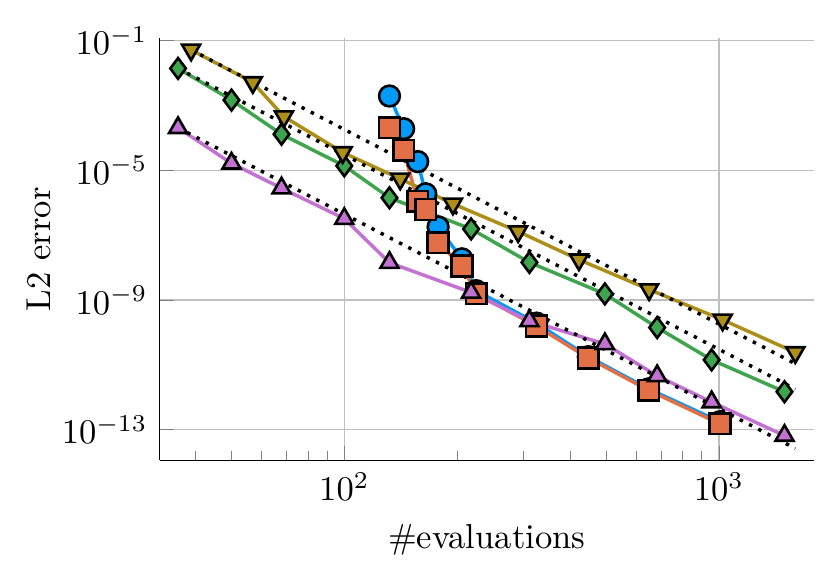}%
    \caption{Logistic Equation}
    \label{fig:dp5:logistic}
  \end{subfigure}
  ~
  \begin{subfigure}{0.48\linewidth}
    \includegraphics[width=\linewidth]{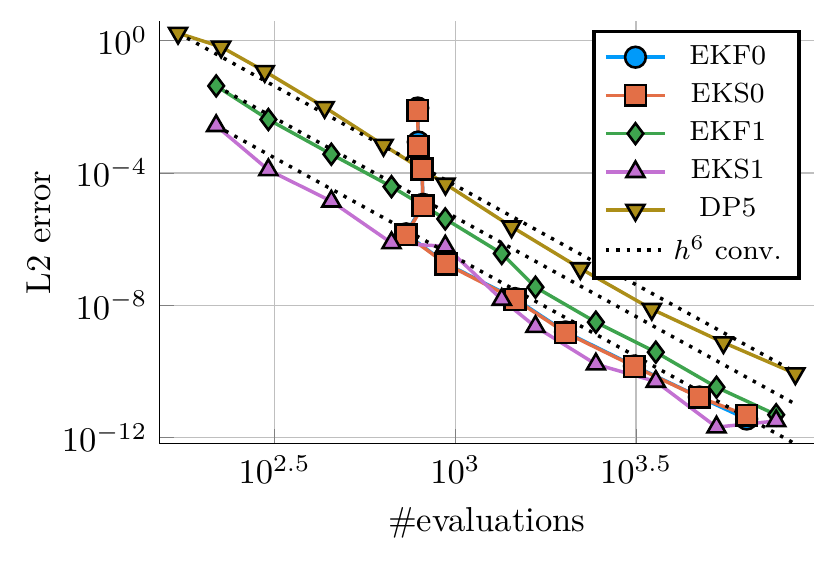}%
    \caption{Lotka-Volterra}
    \label{fig:dp5:lotkavolterra}
  \end{subfigure}
  \\
  \begin{subfigure}{0.48\linewidth}
    \includegraphics[width=\linewidth]{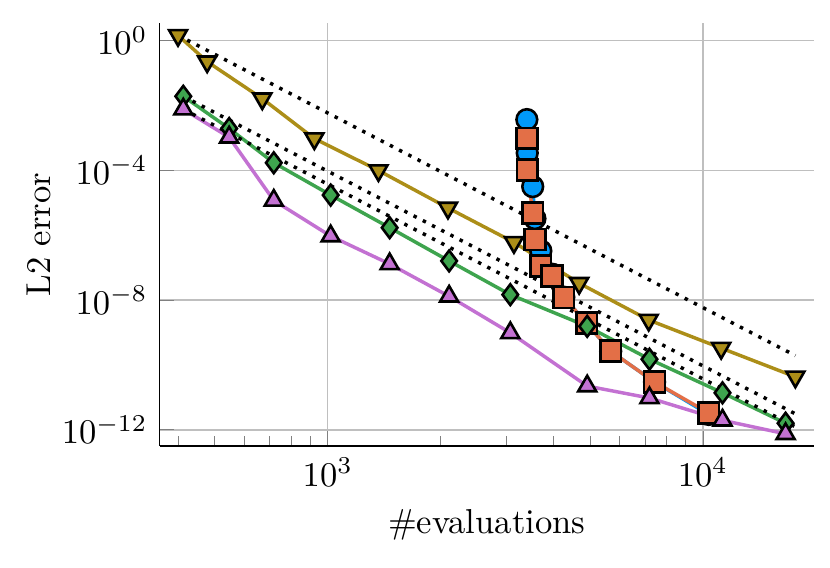}%
    \caption{FitzHugh-Nagumo}
    \label{fig:dp5:fitzhughnagumo}
  \end{subfigure}
  ~
  \begin{subfigure}{0.48\linewidth}
    \includegraphics[width=\linewidth]{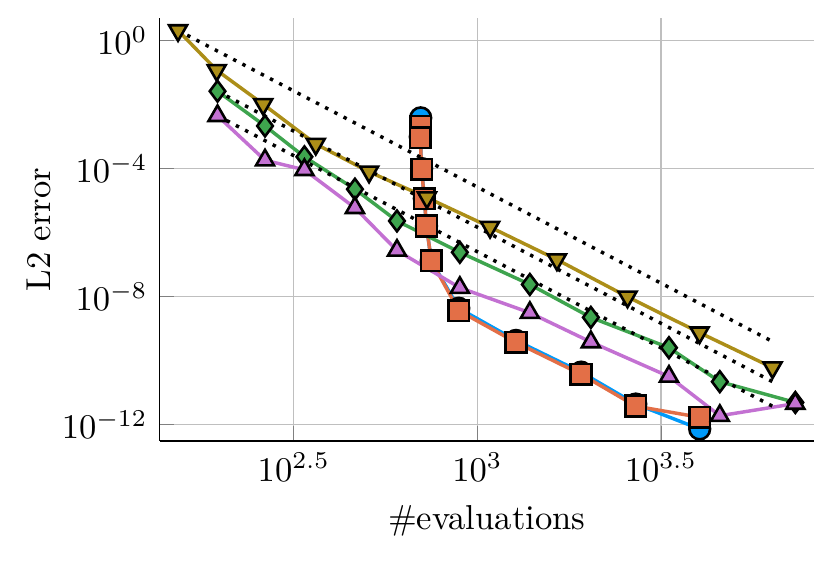}%
    \caption{Brusselator}
    \label{fig:dp5:brusselator}
  \end{subfigure}
  \caption{
    \emph{Comparison to Dormand-Prince 4/5 (DP5)}.
    All probabilistic ODE solvers use an IWP-5 prior and a scalar, time-varying diffusion model.
  }
  \label{fig:dp5}
\end{figure}

Comparison of the probabilistic solvers to the classic Dormand-Prince 4/5 method on additional problems.
\Cref{fig:dp5} presents the resulting work-precision diagrams.
This extends the results of \cref{sec:dp5-comparison}.

\label{sec:supplementary:dp5}

\end{document}